\documentclass[a4paper,10pt,reqno]{amsart}
\usepackage[english]{babel}
\usepackage[centertags]{amsmath}
\usepackage{amsfonts,amssymb,amsthm}
\usepackage{float,longtable}
\usepackage{graphicx}           
\usepackage[latin1]{inputenc}
\usepackage{newlfont}

\newcommand{\hs}{{\mathcal H}_{\sigma}}

\newcommand{\ag}{{\mathcal A}_{\gamma}}
\newcommand{\bg}{{\mathcal B}_{\gamma}}

\newcommand{\R}{{\mathbb R}}

\newcommand{\T}{{\mathbb T}}
\newcommand{\Z}{{\mathbb Z}}
\newcommand{\Q}{{\mathbb Q}}
\newcommand{\Cx}{\mathbb{C}}

\newcommand{\p}{\varphi_1}
\newcommand{\pp}{\varphi_2}
\newcommand{\ph}{\varphi}
\newcommand{\om}{\omega_1}
\newcommand{\omm}{\omega_2}
\newcommand{\norm}[1]{\left\|#1\right\|_{\sigma}}
\newcommand{\normm}[1]{\| #1 \|_{\sigma}}
\newcommand{\e}{\varepsilon}

\newcommand{\lae}{L_{a,\e}}
\newcommand{\ul}{\langle}
\newcommand{\ur}{\rangle}
\newcommand{\av}[1]{\ul #1 \ur}
\newcommand{\avp}[1]{\ul #1 \ur_{\p}}
\newcommand{\avpp}[1]{\ul #1 \ur_{\pp}}
\newcommand{\uz}{\hat{u}_{0,0}}
\newcommand{\str}{\hspace{\stretch{1}}}

\renewcommand{\b}{\beta_0}
\renewcommand{\c}{ \mathrm{cn} }
\newcommand{\am}{ \mathrm{am} }
\newcommand{\x}{\xi}
\renewcommand{\d}{\mathrm{dn}}
\newcommand{\s}{\mathrm{sn}}
\renewcommand{\th}{\vartheta}
\newcommand{\bu}{\bar{u}}
\newcommand{\bv}{\bar{v}}
\newcommand{\bw}{\bar{w}}
\newcommand{\m}{\bar{m}}
\newcommand{\cd}{ \circ}

\newcommand{\mb}{M_{b,\e}}
\newcommand{\mbr}{M_{b,\e |P}}
\newcommand{\V}{\bar{V}}
\newcommand{\Omb}{\bar{\Omega}}

\DeclareMathOperator{\IM}{Im}
\newcommand{\Om}{\Omega}

\newcommand{\jump}{\vspace{9pt}}

\newcommand{\re}[1]{(\ref{#1})}
\begin{document}

\title[Quasi-periodic solutions of $v_{tt} - v_{xx}  +v^3 = f(v)$]
{Quasi-periodic solutions of the equation $v_{tt} - v_{xx} +v^3 = f(v)$}
\author[P. Baldi]{Pietro Baldi}
\address{Pietro Baldi - Sissa,\, via Beirut 2--4,\, 34014 Trieste,\, Italy.}
\email[]{baldi@sissa.it}

\begin{abstract}
We consider 1D completely resonant nonlinear wave equations of the type 
\:$v_{tt} - v_{xx} = -v^3 + \mathcal{O}(v^4)$\: with spatial periodic boundary conditions. We prove the existence of a new type of quasi-periodic small amplitude solutions with two frequencies, for more general nonlinearities. These solutions turn out to be, at the first order, the superposition of a traveling wave and a modulation of long period, depending only on time. 
\end{abstract} 

\maketitle

\section{Introduction}

This paper\footnote{\emph{Keywords:} Nonlinear Wave Equation, Infinite dimensional Hamiltonian Systems, Quasi-periodic solutions, Lyapunov-Schmidt reduction.

\emph{2000AMS Subject Classification:} 35L05, 35B15, 37K50.

Supported by MURST within the PRIN 2004 ``Variational methods and nonlinear differential equations''.}
 deals with a class of one-dimensional completely resonant nonlinear wave equations of the type 
\begin{eqnarray} \label{intro 1}
    \left\{ \begin{array}{l}
    v_{tt} - v_{xx} = -v^3 + f(v)    \\
    v(t,x)=v(t, x+2\pi ), \qquad (t,x)\in  \R ^2 \\
\end{array} \right.
\end{eqnarray}

\noindent 
where $f : \R \rightarrow \R$\, is analytic in a neighborhood of $v=0$\, and $f(v)=\mathcal{O}(v^4)$\, as $v \rightarrow 0$.

In the recent paper \cite{Pro}, M. Procesi  proved the existence of small-amplitude quasi-periodic solutions of \re{intro 1}  of the form
\begin{equation} \label{intro 2}
    v(t,x) = u(\omega_1 t+x , \, \omega_2 t-x),  
\end{equation}

\noindent
where $u$\, is an odd analytic function, $2 \pi$-periodic in both its arguments, and the frequencies \,$\om,\omm \sim 1$\, belong to a Cantor-like set of zero Lebesgue measure. It is assumed that $f$\, is odd and \,$f(v)=\mathcal{O}(v^5)$,\, see Theorem 1 in \cite{Pro}.

These solutions $v(t,x)$\, correspond --- at the first order --- to the superposition of two waves, traveling in opposite directions: 
\[
v(t,x)=\sqrt{\e} \big[ r(\om t+x) +s(\omm t-x)+h.o.t. \big]
\]

\noindent
where 
 \,$\om,\omm =1+\mathcal{O}(\e)$.

Motivated by the previous result, we study in the present paper the existence of quasi-periodic solutions of \re{intro 1} having a different form, namely \begin{equation} \label{intro 3}
    v(t,x) = u(\omega_1 t+x , \, \omega_2 t+x) . 
\end{equation}

\noindent 
Moreover we do not assume $f$\, to be odd.

First of all, we have to consider different frequencies than in \cite{Pro}. Precisely, the appropriate choice for the relationship between the amplitude $\e$\, and the frequencies $\om,\omm$\, turns out to be
\begin{equation}    \label{intro choice}
\om =1+\e +b\e^2, \quad \omm=1+b\e^2 ,
\end{equation}

\noindent
where $b \sim 1/2, \: \  \e \sim 0$.\, This choice leads to look for quasi-periodic solutions $v(t,x)$\, of \re{intro 1} of the form
\begin{equation}  \label{intro 4}
v(t,x)=u\big( \e t, \,(1+b\e^2)t+x \big),
\end{equation}

\noindent
where $(b,\e) \in \R^2$,\:  $\frac{1+b\e^2}{\e} \notin \Q$.\,  
On the contrary, taking in \re{intro 3} frequencies \,$\om=1+\e, \: \: \omm=1+a\e$\, as in \cite{Pro}, no quasi-periodic solutions can be found, see Remark in section 2.
We show that there is no loss of generality passing from \re{intro 3} to \re{intro 4}, because all the possible quasi-periodic solutions of \re{intro 1} of the form \re{intro 3} are of the form \eqref{intro 4}, see Appendix B.

Searching small amplitude quasi-periodic solutions of the form \re{intro 4} by means of the Lyapunov-Schmidt method, leads to the usual system of a range equation and a bifurcation equation. 

The former is solved, in a similar way as \cite{Pro}, by means of the standard Contraction Mapping Theorem, for a set of zero measure of the parameters.
These arguments are carried out in section 4.

In section 5 we study the bifurcation equation, which is infinite-dimen\-sional because we deal with a completely resonant equation. Here new difficulties have to be overcome. 
Since $f$\, is not supposed to be odd, we cannot search odd solutions as in \cite{Pro}, so we look for even solutions. In this way, the bifurcation equation contains a new scalar equation for the average of $u$,\, see [$C$-equation] in \re{M 13 projected problem}, and the other equations contain supplementary terms.

To solve the bifurcation equation we use an ODE analysis; we cannot directly use variational methods as in \cite{BB1},\cite{BB2},\cite{BerPro} because we have to ensure that both components  \,$r,s$\, in \eqref{M 13 projected problem} are non-trivial, in order to prove that the solution $v$\, is actually quasi-periodic.
  
First, we find an explicit solution of the bifurcation equation (Lemma 1) by means of Jacobi elliptic functions (following \cite{BCP},\cite{Pro},\cite{GMP}). 

Next we prove its non-degeneracy (Lemmas 2,3,4); these computations are the heart of the present work. Instead of using a computer assisted proof as in \cite{Pro}, we here employ purely analytic arguments, see also \cite{BCP} (however, our problem 
requires much more involved computations than in \cite{BCP}).
In this way we prove the existence of quasi-periodic solutions of \re{intro 1} of the form \re{intro 4}, see Theorem 1 (end of section 5).

>From the physical point of view, this new class of solutions turns out to be, at the first order, the superposition of a traveling wave (with velocity greater than 1) and a modulation of long period, depending only on time: 
\[
v(t,x)= \e \big[r(\e t)+s((1+b\e^2)t+x) + h.o.t. \big].
\]

Finally, in section 6 we show that our arguments can be also used to extend Procesi result to non-odd nonlinearities, see Theorem 2. 

We also mention that recently existence of quasi-periodic solutions with $n$ frequencies have been proved in \cite{Yuan}. The solutions found in \cite{Yuan} belong to a neighborhood of a solution $u_0(t)$ periodic in time, independent of $x$, so they are different from the ones found in the present paper.

\jump
\textbf{Acknowledgments.} We warmly thank Massimiliano Berti for his daily support,  Michela Procesi and Simone Paleari for some useful discussions.

\jump
\section{The functional setting}
 
We consider nonlinear wave equation \re{intro 1},
\begin{eqnarray*} 
    \left\{ \begin{array}{l}
    v_{tt} - v_{xx} = -v^3 + f(v)    \\
    v(t,x)=v(t, x+2\pi ) 
\end{array} \right.
\end{eqnarray*}
 
\noindent 
where $f$\, is analytic in a neighborhood of $v=0$\, and $f(v)=\mathcal{O}(v^4)$\, as $v \rightarrow 0$.\,  
We look for solutions of the form \re{intro 3},
\begin{equation*} 
    v(t,x) = u(\omega_1 t+x , \, \omega_2 t+x),  
\end{equation*}

\noindent
for $\,( \om, \omm) \in \R^2,$\, $\om, \omm \sim 1$ \: and $u$ \:$2 \pi$-periodic in both its arguments. 
Solutions $v(t,x)$\, of the form \re{intro 3} are \emph{quasi-periodic} in time $t$ when $u$ actually depends on both its arguments and the ratio between the periods is irrational, $\frac{\om}{\omm} \notin \Q$. 

We set the problem in the space $\hs$ defined as follows. Denote $\T=\R/2 \pi \Z \,$ the unitary circle,   $\ph = (\p, \pp) \in \T^2$. If $u$ is doubly $2 \pi$-periodic, $u : \T^2 \rightarrow \R$,\, its Fourier series is 
\begin{equation} \label{M1 fourier}
  u(\ph)=\sum_{(m,n)\in \Z^2} \hat{u}_{mn}\, e^{i m \p} e^{i n \pp}.     
\end{equation}

Let $\sigma >0, \  s \geq 0$. We define $\hs$ as the space of the even $2 \pi$-periodic functions $u : \T^2 \rightarrow \R$\, which satisfy 
\begin{equation*} \label{M2 def di hs}
\sum_{(m,n)\in \Z^2} \left|\hat{u}_{mn}\right|^2  \left[ 1 + (m^2 + n^2)^s \right] \, e^{2 \sqrt{m^2 + n^2} \, \sigma} := \left\|u \right\|_{\sigma}^2 < \infty.     
\end{equation*}

The elements of $\hs$ are even periodic functions which admit an analytic extension to the complex strip $\{z \in \Cx : \left| \IM (z) \right| < \sigma \}$.

$(\hs , \left\| \cdot \right\|_{\sigma})$\, is a Hilbert space; for $s > 1\,$ it is also an algebra, that is, there exists a constant $c>0$\, such that  
\[
\norm{uv}\leq c \norm{u} \norm{v} \quad \forall \, u,v \in \hs,
\]

\noindent
see Appendix A. Moreover the inclusion \,$\mathcal{H}_{\sigma,s+1} \hookrightarrow \mathcal{H}_{\sigma,s}$\, is compact. 

We fix $s > 1$ once and for all.

\jump
We note that all the possible quasi-periodic solutions of \re{intro 1} of the form \re{intro 3} are of the form \eqref{intro 4} if we choose frequencies as in \re{intro choice}, see Appendix B. So we can look for solutions of \re{intro 1} of the form \re{intro 4}, \,$v(t,x)=u\big( \e t, \,(1+b\e^2)t+x \big)$,\, without loss of generality. 
For functions of the form \re{intro 4}, problem \re{intro 1} is written as
\begin{eqnarray*} \label{M3}
  \left\{ \! \begin{array}{l}
\e \left[ \e \,\partial_{\p}^2  + 2(1+b\e^2)\, \partial_{\p \pp}^2  + b\e (2+b\e^2)\, \partial_{\pp}^2 \right] (u) =  -u^3+f(u) \vspace{3pt} \\
    u \in \hs .\\
\end{array} \right.
\end{eqnarray*}

We define \,$M_{b,\e}=\e \,\partial_{\p}^2  + 2(1+b\e^2)\, \partial_{\p \pp}^2  + b\e (2+b\e^2)\, \partial_{\pp}^2$,\: rescale \,$u \rightarrow \e u$\, and set 
\,$f_\e(u)=\e^{-3} f(\e u)$,\, so \re{intro 1} can be written as
\begin{eqnarray} \label{M4}
  \left\{ \! \begin{array}{l}
M_{b,\e}[u] =  -\e u^3+ \e f_\e (u) \vspace{3pt} \\
    u \in \hs . \\
\end{array} \right.
\end{eqnarray}

The main result of the present paper is the existence of solutions $u_{(b,\e)}$ of \re{M4} for $(b,\e)$ in a suitable uncountable set (Theorem 1).

\jump
\noindent
\textbf{\emph{Remark.}} 
If we simply choose frequencies \,$\om=1+\e, \: \: \omm=1+a\e$\, as in \cite{Pro}, we obtain a bifurcation equation different than \re{M 13 projected problem}. Precisely, it appears 0\, instead of \,$-b(2+b\e^2)\, s''$\, in the left-hand term of the $Q_2$-equation in \re{M 13 projected problem}; so we do not find solutions which are non-trivial in both its arguments, but only solutions depending on the variable $\p$. This is a problem because the quasi-periodicity condition requires dependence on both variables. 

So we have to choose frequencies depending on $\e$\, in a more general way; a good  choice is \re{intro choice}, \:$\om=1+\e+b\e^2$,\:\,  $\omm= 1+b\e^2 $.

\jump
\section{Lyapunov-Schmidt reduction}

The operator $M_{b,\e}$ is diagonal in the Fourier basis \,$e_{mn} = e^{im\p}e^{in\pp}$\, with eigenvalues $-D_{b,\e}(m,n)$,\, that is, if $u$\, is written in Fourier series as in (\ref{M1 fourier}),
\begin{eqnarray} \label{M5}
	M_{b,\e}[u]= \,- \hspace{-6pt} \sum_{(m,n)\in \Z ^2}D_{b,\e}(m,n)\, \hat{u}_{mn}\, e^{i m \p}\, e^{i n \pp},
\end{eqnarray}

\noindent
where the eigenvalues $D_{b,\e}(m,n)$ are given by
\begin{align} \label{M6}
  D_{b,\e}(m,n) = & \: \e\, m^2 + 2(1+b\e^2)\, mn + b\e (2+b\e^2) \, n^2   \nonumber \\
                = & \: (2+b\e^2)\,\Big( \frac{\e}{2+b\e^2}\,m +n \Big) \big( m+ b\e\,n \big). 
\end{align}

For $\e=0 \ $ the operator is $M_{b,0}=2 \,\partial_{\p \pp}^2$; its kernel $Z$\, is the subspace of functions of the form \,$u(\p,\pp)=r(\p)+s(\pp)$\, for some \,$r,\,s \in \hs$\, one-variable functions, 
\[
Z=\big\{u\in \hs : \hat{u}_{mn}=0 \ \  \forall \, (m,n) \in \Z^2, \,m,n \neq 0 \big\}.
\]

We can decompose $\hs$ in four subspaces setting
\begin{eqnarray} \label{M8 decomp}
 \begin{array}{cl}	
 &	\quad C=\{u \in \hs : u(\ph)=\hat{u}_{0,0} \} \cong \R, \vspace{2pt}\\
 &	\quad Q_1 = \{u \in \hs : u(\ph)= \sum_{m \neq 0} \hat{u}_{m,0}\,e^{im\p} =r(\p) \}, \vspace{2pt}\\
 &	\quad Q_2 = \{u \in \hs : u(\ph)= \sum_{n \neq 0} \hat{u}_{0,n}\,e^{in\pp} =s(\pp) \}, \vspace{2pt}\\
 &	\quad P = \{u \in \hs : u(\ph)= \sum_{m,n \neq 0} \hat{u}_{mn}\,e^{im\p}\,e^{in\pp} =p(\p, \pp) \} .
 \end{array}
\end{eqnarray}

\noindent
Thus the kernel is the direct sum \,$Z=C\oplus Q_1 \oplus Q_2 \,$ and the whole space is $\hs = Z \oplus P$. Any element $u$\, can be decomposed as 
\begin{eqnarray} \label{M9 decomp2}
  \begin{array}{rl}
	u(\ph)= & \hspace{-4pt} \hat{u}_{0,0} + r(\p) + s(\pp) + p(\p,\pp) \\
	 = & \hspace{-4pt} z(\ph)+p(\ph).
  \end{array}
\end{eqnarray}
  
We denote $\langle \: \cdot \: \rangle$ the integral average: given $g \in \hs$, 
\begin{eqnarray*} \label{M 10 average}
  \begin{array}{c}
    \langle g \rangle =      \ul g \ur_{\ph} = \frac{1}{(2 \pi)^2} \int_0^{2 \pi} \int_0^{2 \pi} g(\ph) \, d\p d \pp ,\vspace{7pt}\\
    \ul g \ur_{\p} =  \frac{1}{2 \pi} \int_0^{2 \pi} g(\ph) \, d \p, \qquad 
    \ul g \ur_{\pp} =  \frac{1}{2 \pi} \int_0^{2 \pi} g(\ph) \, d \pp .
  \end{array}  
\end{eqnarray*}   

\noindent    
Note that $\frac{1}{2 \pi} \int_0^{2 \pi} e^{ikt} \, d t =0 \,$ for all integers $k \neq 0, \,$ so 
\begin{eqnarray*} \label{M 11 averages}
    \begin{array}{lcl} 
      \av{r}=\avp{r}=0 & & \avpp{r}=r \vspace{2pt}\\    
      \av{s}=\avpp{s}=0 & & \avp{s}=s \vspace{2pt}\\
      \av{p}=\avp{p}=\avpp{p}=0 & & \av{u}=\hat{u}_{0,0} \vspace{2pt}\\         
    \end{array}
\end{eqnarray*} 

\noindent
for all $r \in Q_1, \, s \in Q_2, \, p \in P, \, u \in \hs, \,$ and by means of these averages we can construct the projections on the subspaces,
\begin{eqnarray*} \label{M 12 proj}
	 %\hspace{30pt}
	   \Pi_{C}=\av{\, \cdot \,}, \qquad
	   \Pi_{Q_1}=\avpp{\, \cdot \,} - \av{\, \cdot \,}, \qquad
	   \Pi_{Q_2}=\avp{\, \cdot \,} - \av{\, \cdot \,}. 
\end{eqnarray*}

Let $u=z+p \,$ as in  \re{M9 decomp2}; we write $u^3\,$ as $u^3=z^3+(u^3-z^3)\,$ and compute the cube $z^3=(\hat{u}_{0,0}+r+s)^3$. The operator $\mb$\, maps every subspace of \re{M8 decomp} in itself and it holds \,$ \mb [r]=\e r''$,\,  
$\mb [s]= b\e(2+b\e^2)s''$,\, 
$\mb [\hat{u}_{0,0}]=0$.\, So we can project our problem \re{M4} on the four  subspaces:
\begin{eqnarray} \label{M 13 projected problem}
      \begin{array}{rl}
       0= & \hspace{-4pt} \uz ^3+3\uz\left(\av{r^2}+\av{s^2}\right)+\av{r^3}
                                                       +\av{s^3}+ \vspace{5pt}\\
          &+\Pi_C\left[(u^3-z^3)-f_{\e}(u)\right]\vspace{17pt}
                                \str  \left[ C\text{\emph{-equation}} \right]\\
  -r''= & \hspace{-4pt}3\uz^2 r+3\uz \left(r^2 - \av{r^2}\right)+r^3 -\av{r^3}
                                                      +3\av{s^2}\,r+\vspace{5pt}\\            &     +\Pi_{Q_1}\left[(u^3-z^3)-f_{\e}(u)\right] \vspace{17pt}
                                  \str  \left[ Q_1\text{\emph{-equation}} \right]\\
  -b(2+b\e^2)\, s'' = & \hspace{-4pt} 3\uz^2 s+3\uz \left( s^2 - \av{s^2} \right) 
                                        +s^3 -\av{s^3} + 3\av{r^2}\,s+\vspace{5pt}\\ 
          &    +\Pi_{Q_2}\left[(u^3-z^3)-f_{\e}(u)\right] \vspace{17pt}
                                 \str  \left[ Q_2\text{\emph{-equation}} \right]\\
  M_{b,\e} [p]= & \hspace{-4pt} \e \, \Pi_P \left[-u^3 +f_{\e}(u)\right].\vspace{5pt}
                           \str  \left[ P \text{\emph{-equation}} \right]\vspace{5pt}\\
  \end{array}   
\end{eqnarray}

Now we study separately the projected equations.

\jump
\section{The range equation}

We write the $P$-equation thinking $p\,$ as variable and $z\,$ as a ``parameter'', 
\begin{eqnarray*} \label{M 14 P equation}
	\mb [p]=  \e \, \Pi_P \left[-(z+p)^3 +f_{\e}(z+p)\right].
\end{eqnarray*}

We would like to invert the operator $\mb$. In Appendix C we prove that, fixed any $\gamma \in (0,\frac{1}{4}), \,$ there exists a non-empty uncountable set $\bg \subseteq \R^2 \,$ such that, for all $(b,\e) \in \bg, \,$ it holds
\begin{eqnarray*} \label{M 15 stima per D_be}
 \left| D_{b,\e}(m,n) \right| >   \gamma \qquad \forall \, m,n \in \Z,\: m,n \neq 0.	
\end{eqnarray*}

\noindent
Precisely, our Cantor set $\bg$\, is 
\begin{equation*}  \label{vero bg nel testo}
\bg=\bigg\{ (b,\e) \in \R^2 : \: \: \frac{\e}{2+b\e^2},b\e^2 \in \tilde{B}_\gamma,\ \ 
\Big| \frac{\e}{2+b\e^2} \Big|, |b\e^2 | < \frac{1}{4}, \ \ \frac{1+b\e^2}{\e} \notin \Q \bigg\}, 
\end{equation*}

\noindent
where $\tilde{B}_\gamma$\, is a set of ``badly approximable numbers'' defined as
\begin{eqnarray} \label{def di B_gamma nel testo}
	\tilde{B}_\gamma= \Big\{ x\in \R : \left| m+nx \right| > \frac{\gamma}{| n | } \ \ \forall \, m,n \in \Z,\, m\neq 0,\, n \neq 0 \Big\},
\end{eqnarray}

\noindent
see Appendix C. Therefore $M_{b,\e|P}\,$ is invertible for $(b,\e) \in \bg \,$ and by \re{M5} it follows
\begin{eqnarray*} \label{M 16}
   ( \mbr ) ^{-1}[h]= \,- \hspace{-4pt} \sum_{m,n \neq 0} \frac{\hat{h}_{mn}}{D_{b,\e}(m,n)} \,e^{im\p}\, e^{in\pp} 
\end{eqnarray*}

\noindent
for every $h=\sum_{m,n \neq 0}\hat{h}_{mn} \,e^{im\p}\, e^{in\pp} \in P$. Thus we obtain a bound for the inverse operators, uniformely in $(b,\e) \in \bg$: 
\begin{eqnarray*} \label{M 17}
	\left\| ( \mbr )^{-1} \right\| \leq \frac{1}{\gamma}.
\end{eqnarray*}

Applying the inverse operator $( \mbr )^{-1}$,\, the $P$-equation becomes 
\begin{eqnarray} \label{M 18 inverse P equation}
	p+\e ( \mbr )^{-1} \Pi_P \left[(z+p)^3 -f_{\e}(z+p)\right]=0.
\end{eqnarray}

We would like to apply the Implicit Function Theorem, but the inverse operator $( \mbr )^{-1}\,$ is defined only for $(b,\e) \in \bg \,$ and in the set $\bg\,$ there are infinitely many holes, see Appendix C. So we fix $(b,\e) \in \bg, \,$  introduce an auxiliary parameter $\mu\,$ and consider the auxiliary equation
\begin{eqnarray} \label{M 19 aux P eq}
 p+\mu ( \mbr )^{-1} \Pi_P \left[(z+p)^3 -f_{\mu}(z+p)\right]=0.
\end{eqnarray}

Following Lemma 2.2 in \cite{Pro}, we can prove, by the standard Contraction Mapping Theorem, that there exists a positive constant $c_1$\, depending only on $f$ such that, if 
\begin{equation}  \label{dominio di p implicita}
(\mu,z) \in \R \times Z, \qquad |\mu| \norm{z}^2 < c_1 \gamma,
\end{equation}

\noindent
equation \re{M 19 aux P eq} admits a solution $p_{(b,\e)}(\mu,z) \in P$. 
Moreover, there exists a positive constant $c_2$\, such that the solution $p_{(b,\e)}(\mu,z)$\, respects the bound
\begin{equation}  \label{bound for p}
\norm{p_{(b,\e)}(\mu,z)} \leq \frac{c_2}{\gamma} \norm{z}^3 |\mu|.
\end{equation} 

\noindent
Than we can apply the Implicit Function Theorem to the operator
\begin{align*} \label{M 20 aux F(mu,z,p)}
& \R \times Z \times P \longrightarrow P  \\
& (\mu,z,p) \longmapsto p+\mu ( \mbr )^{-1} \Pi_P \left[(z+p)^3 -f_{\mu}(z+p)\right]
\end{align*}

\noindent
at every point $(0,z,0)$, so, by local uniqueness, we obtain the regularity: $p_{(b,\e)}$,\, as function of $(\mu,z)$,\, is at least of class $\mathcal{C}^1$. 

Notice that the domain of any function $p_{(b,\e)}$\, is defined by \re{dominio di p implicita}, so it does not depend on $(b,\e) \in \bg$. 
   
In order to solve \re{M 18 inverse P equation}, we will need to evaluate  $p_{(b,\e)}$\, at $\mu=\e$; we will do it as last step, after the study of the bifurcation equation. 

We observe that in these computations we have used the Hilbert algebra property of the space $\hs, \ \, \norm{uv} \leq c \norm{u}\norm{v}\, \forall \, u,v \in \hs.$

\jump
\section{The bifurcation equation}
We consider auxiliary $Z$-equations: we put $f_{\mu}\,$ instead of $f_{\e}\,$ in \re{M 13 projected problem},

\begin{eqnarray} \label{M 22 aux Z eq}
   \begin{array}{rl}
       0= & \hspace{-4pt} \uz ^3+3\uz\left(\av{r^2}+\av{s^2}\right)+\av{r^3}
                                                       +\av{s^3}+ \vspace{5pt}\\
          &+\Pi_C\left[(u^3-z^3)-f_{\mu}(u)\right]\vspace{17pt}
                                \str  \left[ C-equation \right]\\
  -r''= & \hspace{-4pt}3\uz^2 r+3\uz \left(r^2 - \av{r^2}\right)+r^3 -\av{r^3}
                                                      +3\av{s^2}\,r+\vspace{5pt}\\            &     +\Pi_{Q_1}\left[(u^3-z^3)-f_{\mu}(u)\right] \vspace{17pt}
                                  \str  \left[ Q_1-equation \right]\\
  -b(2+b\e^2)\, s'' = & \hspace{-4pt} 3\uz^2 s+3\uz \left( s^2 - \av{s^2} \right) 
                                        +s^3 -\av{s^3} + 3\av{r^2}\,s+\vspace{5pt}\\ 
          &    +\Pi_{Q_2}\left[(u^3-z^3)-f_{\mu}(u)\right]. \vspace{5pt}
                                 \str  \left[ Q_2 -equation \right]\\
  \end{array}   
\end{eqnarray}

We substitute the solution $p_{(b,\e)}(\mu,z)\,$ of the auxiliary $P$-equation \re{M 19 aux P eq} inside the auxiliary $Z$-equations \re{M 22 aux Z eq}, writing \,$u=z+p=z+p_{(b,\e)}(\mu,z)$,\, for $(\mu,z)$\, in the domain \re{dominio di p implicita} of $p_{(b,\e)}$.

We have $\, p_{(b,\e)}(\mu,z) \!=\! 0\ $ for $\mu\!=\!0, \  $ so the term $\left[(u^3-z^3)-f_{\mu}(u)\right]$ vanishes for $\mu=0$\, and the bifurcation equations at $\mu=0 \,$ become
\begin{eqnarray} \label{M 23 bif eq for mu=0}
  \begin{array}{lrl}
%\hspace{10pt} 
&      0= & \hspace{-4pt}\uz^3+3\uz\left(\av{r^2}+\av{s^2}\right)+
                          \av{r^3}+\av{s^3} \vspace{6pt} \hspace{10pt}\\
%\hspace{10pt} 
& -r''= & \hspace{-4pt}3\uz^2 r+3\uz \left(r^2 -                                   \av{r^2}\right)+r^3 -\av{r^3}+3\av{s^2}\,r \vspace{6pt} \hspace{10pt}\\  
%\hspace{10pt} 
& -b(2+b\e^2)\, s'' = & \hspace{-4pt} 3\uz^2 s+3\uz \left( s^2 - \av{s^2}                                 \right) +s^3 -\av{s^3} + 3\av{r^2}\,s . \hspace{10pt}
  \end{array}    
\end{eqnarray} 
\vspace{0pt}

We look for non-trivial $z=\uz+r(\p)+s(\pp) \,$ solution of \re{M 23 bif eq for mu=0}. We  rescale setting 
\begin{eqnarray} \label{M 24 magic rescaling}  
	\begin{array}{ll}
	r=  x  & \hspace{30pt} \uz= c \\
	s= \sqrt{b(2+b\e^2)}\ y  & \hspace{30pt} \lambda=\lambda_{b,\e}=b(2+b\e^2),  \\
	\end{array}
\end{eqnarray}

\noindent 
so the equations become \vspace{3pt}
\begin{eqnarray}    	\label{M 25 il cuore}
 \begin{array}{ll}
 \quad	&  c^3+3c \big( \av{x^2}+\lambda \av{y^2} \big)        +\av{x^3}+\lambda^{3/2}\av{y^3}=0 \vspace{6pt}\\
 \quad	&  x''+3c^2 x+3c \big(x^2 - \av{x^2} \big)+x^3-\av{x^3}+3\lambda \av{y^2}\,x =0  \vspace{6pt}\\
 \quad	&  y''+3c^2 \frac{1}{\lambda}y+3c\frac{1}{\sqrt{\lambda}} \big(y^2 - \av{y^2} \big)+y^3-\av{y^3}+3\frac{1}{\lambda} \av{x^2}\,y =0 .
	\end{array}  
\end{eqnarray}

In the following we show that, for $|\lambda -1|\,$ sufficiently small, the system \re{M 25 il cuore} admits a non-trivial non-degenerate solution. We consider $\lambda\,$ as a free real parameter, recall that $\,Z=C \times Q_1 \times Q_2$\:  and define 
 \,$  G:\R \times Z \rightarrow Z$\:  
setting \,$G(\lambda, c,x,y)\,$ as the set of three left-hand terms of \re{M 25 il cuore}.

\jump
\noindent \textbf{Lemma 1.} \emph{There exist $\bar{\sigma}>0$\, and a non-trivial one-variable even analytic function $\b$\, belonging to $\hs$\, for every $\sigma \in (0,\bar{\sigma})$, such that \,$G(1,0,\b,\b)=0$,\: that is \,$(0,\b,\b)$\, solves \re{M 25 il cuore} for $\lambda=1$.}

\jump
\noindent
\emph{Proof.} We prove the existence of a non-trivial even analytic function $\b$\, which satisfies
\begin{eqnarray} \label{M 26 eq di beta 0}
	\b''+\b^3+3\av{\b^2}\, \b=0, \qquad   \av{\b}=\av{\b ^3}=0.
\end{eqnarray}

For any  
$m\in(0,1)$\, we consider the Jacobi amplitude $\am(\,\cdot\,,m):\R \rightarrow \R\,$ as the inverse of the elliptic integral of the first kind 
\[
I(\,\cdot\,,m):\R \rightarrow \R, \qquad I(\ph,m)=\int_0^{\ph} \frac{d\vartheta}{\sqrt{1-m \sin^2 \vartheta}}.
\]

\noindent
We define the Jacobi elliptic cosine setting 
\[
\c(\xi)=\c(\x,m)=\cos(\am(\xi,m)),
\]

\noindent
see \cite{Handbook} ch.16, \cite{MathWorld}. It is a periodic function of period $4K$,\, where $K=K(m)\,$ is the complete elliptic integral of the first kind
\[
K(m)=\int_0^{\pi /2} \frac{d\vartheta}{\sqrt{1-m \sin^2 \vartheta}}.
\]

Jacobi cosine is even, and it is also odd-symmetric with respect to $K\,$ on $[0, 2K],\,$ that is $\,\c (\xi +K)=-\c (\xi -K)$,\,  just like the usual cosine. Then the averages on the period $4K$\, are 
\begin{eqnarray*} \label{App2 averages}
	\av{\c}=\av{\c^3}=0.
\end{eqnarray*}

 \noindent
Therefore it admits an analytic extension with a pole at $iK'$,\, where $K'=K(1-m)$,\, and it satisfies 
$\,	(\c')^2=-m \,\c^4 +(2m-1)\, \c^2 + (1-m),\ $ then $\c$\, is a solution of the ODE
\begin{eqnarray*} \label{cn eq}
	\c''+2m \,\c^3 +(1-2m)\,\c =0.
\end{eqnarray*}

We set $\,\b(\x)=V\c(\Om \x ,m)\:$ for some real parameters $V,\Om > 0, \, m \in (0,1)$.\, $\b$\, has a pole at $i\frac{K'}{\Omega}$,\, so it belongs to $\hs$\, for every $0<\sigma<\frac{K'}{\Omega}$. 

$\b$\, satisfies
\begin{eqnarray*} \label{App2 1}
	\b''+\Big( 2m \frac{\Om^2}{V^2}\Big)\, \b^3 + \Om^2 (1-2m)\,\b=0.
\end{eqnarray*} 

\noindent
If there holds the equality $\, 2m\Om^2=V^2,\,$ the equation becomes
\begin{eqnarray*} \label{App2 2}
	\b''+\b^3 + \Om^2 (1-2m)\,\b=0.
\end{eqnarray*}

$\b$\, is \,$\frac{4K(m)}{\Om}$-periodic; it is \,$2\pi$-periodic if \,$\Om=\frac{2K(m)}{\pi}$.\, Hence we require
\begin{eqnarray} \label{App2 rel}
	2m\Om^2=V^2, \qquad \Om=\frac{2K(m)}{\pi}.
\end{eqnarray}

The other Jacobi elliptic functions we will use are 
\[
\mathrm{sn}(\x)=\sin (\am(\x,m)), \qquad \mathrm{dn}(\x)=\sqrt{1-m \,\mathrm{sn}^2(\x)},
\]

\noindent
see \cite{Handbook},\cite{MathWorld}. From the equality \,$m\, \c^2(\x)=\mathrm{dn}^2(\x)-(1-m),$\, with change of variable \,$x=\am(\x)$\, we obtain 
\[
\int_0^{K(m)} m\,\c^2(\x)\, d\x =E(m)-(1-m)K(m), 
\]

\noindent
where $E(m)\,$ is the complete elliptic integral of the second kind, 
\[
E(m)=\int_0^{\pi /2} \sqrt{1-m \sin^2 \vartheta}\, d\vartheta.
\]

\noindent
Thus the average on $[0,2 \pi]$\, of $\b^2$\, is
\begin{eqnarray*}
	\av{\b^2}=\frac{V^2}{m \,K(m)}\big[E(m)-(1-m)K(m)\big].
\end{eqnarray*}

\noindent
We want the equality $\,3\av{\b^2}=\Om^2 (1-2m)$\, and this is true if
\begin{eqnarray} \label{App2 numerico...}
	E(m)+\frac{8m-7}{6}\, K(m)=0.
\end{eqnarray}

\noindent The left-hand term \,$\psi(m):=E(m)+\frac{8m-7}{6}\, K(m)$\, is continuous in $m$; its value at $m=0$\, is \,$-(\pi/12) <0$,\, while at $m=1/2$,\, by definition of $E$\, and $K$,
\[
\psi\Big(\frac{1}{2}\Big)=\frac{1}{2}\int_0^{\pi/2} \frac{\cos^2 \th}{ \big(1-\frac{1}{2}\sin^2 \th \big)^{1/2}  }\, d\th >0.
\]

\noindent
Moreover, its derivative is strictly positive for every $m \in \big[0,\frac{1}{2}\big]$,
\[
\begin{array}{rl}
\psi'(m)= & \hspace{-4pt}  \int_0^{\pi/2}\, \frac{ 8-\frac{5}{2}\sin^2 \th +3m \sin^4 \th -8m \sin^2 \th }{ 6(1-m\sin^2 \th)^{3/2}    } \, d\th \vspace{6pt}\\
\geq & \hspace{-4pt}\int_0^{ \pi/2 }\, \frac{3+\cos^2 \th}{6} \, d\th \, >0,
\end{array}
\]

\noindent
hence there exists a unique $\bar{m} \in \big(0,\frac{1}{2}\big)$\, which solves \re{App2 numerico...}. Thanks to the tables in \cite{Handbook}, p. 608-609, we have $0.20<\bar{m}<0.21$.

By \re{App2 rel} the value $\bar{m}$\, determines the parameters $\,\bar{\Om}$\, and $\,\bar{V}$,\, so the function $\b(\x)= \bar{V} \c (\bar{\Om}\x, \bar{m})$\, satisfies \re{M 26 eq di beta 0} and 
 \,$(0,\b,\b)$\, is a solution of \re{M 25 il cuore} for $\lambda=1$.  
Therefore $\b \in \hs$\, for every $\sigma \in (0,\bar{\sigma})$,\, where 
$ \bar{\sigma}=( \frac{K'}{\Omega} )_{|m=\bar{m}}. $ 
\str $\Box$

\jump
The next step will be to prove the non-degeneracy of the solution $(1,0,\b,\b)$,\, that is to show that the partial derivative \,$\partial_{Z}G(1,0,\b,\b)$\, is an invertible operator. 
This is the heart of the present paper. 
We need some preliminary results. 

\jump 
\noindent \textbf{Lemma 2.} \emph{Given $h$\, even $2\pi$-periodic, there exists a unique even $2\pi$-periodic  $w$\, such that} 
\[
w''+\big( 3\b^2 +3\av{\b^2} \big) w =h.
\]

\noindent
\emph{This defines the Green operator \,$L: \hs \rightarrow \hs$, \:\:$L[h]=w$.} 

\jump
\noindent
\emph{Proof.} 
We fix a $2\pi$-periodic even function $h$.\, We look for even $2\pi$-periodic solutions of the non-homogeneous equation
\begin{eqnarray} \label{App3 1}
	x''+\big( 3\b^2 +3\av{\b^2} \big) x =h.
\end{eqnarray}

First of all, we construct two solutions of the homogeneous equation
\begin{eqnarray} \label{App3 hom}
	x''+\big( 3\b^2 +3\av{\b^2} \big) x =0.
\end{eqnarray}

We recall that $\b$\, satisfies $\,\b''+\b^3+3\av{\b^2}\, \b=0$,\, then deriving with respect to its argument $\x$\, we obtain $\,\b'''+3\b^2 \b' +3\av{\b^2}\b'=0,$\, so $\b'$\, satisfies \re{App3 hom}. We set 
\begin{eqnarray} \label{App3 2}
	\bar{u}(\x)=-\frac{1}{\bar{V} \bar{\Om}^2}\b'(\x)=-\frac{1}{\bar{\Om}}\, \c '(\bar{\Om} \x,\bar{m}),
\end{eqnarray}

\noindent
thus $\bar{u}$\, is the solution of the homogeneous equation such that $\bar{u}(0)=0$, \ $\bar{u}'(0)=1$.\, It is odd and $2\pi$-periodic.

Now we construct the other solution. We indicate $c_0$\, the constant $c_0 =\av{\b}$.\, We recall that, for any $V,\Om,m$\, the function $y(\x)=V\c (\Om \x ,m)$\, satisfies  
\begin{eqnarray*} \label{App3 nuova 1}
	y''+\Big( 2m \frac{\Om^2}{V^2}\Big)\, y^3 + \Om^2 (1-2m)\,y=0.
\end{eqnarray*} 

\noindent
We consider $m$\, and $V$\, as functions of the parameter $\Om$,\, setting
\begin{equation}   \label{App3 nuova 2, rel}
m=m(\Om)=\frac{1}{2} - \frac{3 c_0}{2 \Om^2}, \qquad V=V(\Om)= \sqrt{\Om^2-3c_0}.
\end{equation}  

\noindent
We indicate \,$y_\Om (\x)=V(\Om)\, \c \big( \Om \x, m(\Om) \big),$\, 
so \,$(y_\Om)_\Om$\, is a one-parameter family of solutions of
\begin{equation*}  \label{App3 nuova 4}
y_\Om '' + y_\Om ^3 + 3c_0 y_\Om =0.
\end{equation*}

\noindent
We can derive this equation with respect to $\Om$,\, obtaining
\begin{equation*}  \label{App3 nuova 5}
(\partial_\Om y_\Om )''+3 y_\Om ^2 (\partial_\Om y_\Om )+3c_0 (\partial_\Om y_\Om )=0.
\end{equation*}

\noindent
Now we evaluate \,$(\partial_\Om y_\Om )$\, at \,$\Om=\bar{\Om}$,\, where $\bar{\Om}$\, correspond to the value $\bar{m}$\, found in Lemma 1. For $\Om=\bar{\Om}$\, it holds  
\,$y_{\bar{\Om}}=\b$, so \,$(\partial_\Om y_\Om )_{|\Om=\bar{\Om}}$\, satisfy \re{App3 hom}. In order to normalize this solution, we compute
\[
(\partial_\Om y_\Om )(\x)= (\partial_\Om V )\c(\Om \x, m) + V \x \c'(\Om \x, m)+ V  \partial_m \c(\Om \x,m)(\partial_\Om m).
\]

\noindent
Since $\c(0,m)=1$\, $\forall \,m$,\, it holds \,$\partial_m \c(0,m)=0$;\, therefore $\c'(0,m)=0$\, $\forall \,m$. From \re{App3 nuova 2, rel} we have \,$\partial_\Om V =\frac{\Om}{V}$,\, so we can normalize setting
\[
\bar{v}(\x)=\frac{\bar{V}}{\bar{\Om}}\,\big(\partial_\Om y_\Om \big)_{|\Om=\bar{\Om}} (\x).
\]    

%\noindent
$\bar{v}$\, is the solution of the homogeneous equation \re{App3 hom} such that $\bar{v}(0)=1$,\, $\bar{v}'(0)=0$. 
We can write an explicit formula for $\bar{v}$. From the definitions it follows for any $m$
\[
\partial_m \am (\x,m)=-\d (\x,m) \frac{1}{2} \int_0^{\x} \frac{\s^2(t,m)}{\d^2(t,m)}\,dt.
\]

\noindent
Therefore $\c'(\x)=-\s(\x) \,\d(\x)$;\, then we obtain for $(V,\Om,m)=(\bar{V},\bar{\Om},\bar{m})$ 
\begin{eqnarray} \label{App3 formula magica}
\bar{v}(\x)= \c(\bar{\Om} \x) +\frac{\bar{V}^2}{\bar{\Om}} \, \c'(\bar{\Om} \x) \, \bigg[ \,\x +\, \frac{2\bar{m}-1}{2}\, \int_0^{\x} \frac{\s^2(\bar{\Om} t)}{\d^2(\bar{\Om} t)}\,dt \bigg].
 \end{eqnarray}

By formula \re{App3 formula magica} we can see that $\bv$\, is even; it is not periodic and there holds
\begin{eqnarray} \label{App3 diff}
	\bv(\x+2\pi)-\bv(\x)=\frac{\bar{V}^2 k}{\bar{\Om}}\, \c'(\bar{\Om} \x)= \,-\bar{V}^2 k \,\bu (\x) ,  
\end{eqnarray}

\noindent
where 
\begin{equation} \label{App3 def di k}
k:= \,2\pi+\, \frac{2\bar{m}-1}{2}\, \int_0^{2\pi} \frac{\s^2(\bar{\Om} t)}{\d^2(\bar{\Om} t)}\,dt.
\end{equation}

\noindent
>From the equalities (L.1) and (L.2) of Lemma 3\,  we obtain 
\begin{equation} \label{App3 polinomio per k}
k=2\pi\, \frac{-1+16\bar{m}-16\bar{m}^2}{12\bar{m}(1-\bar{m})} ,
\end{equation}

\noindent
so $k>0$\, because $\bar{m} \in (0.20, 0.21)$.

\jump 
We have constructed two solutions $\, \bu, \, \bv \,$ of the homogeneous equation; their wronskian \,$\bu' \bv -\bu \bv'$\, is equal to $1$,\, so we can write a particular solution $\bw$\, of the non-homogeneous equation \re{App3 1} as
\begin{eqnarray*} \label{App3 bar w}
	\bw(\x)=\bigg( \int_0^{\x}h \bv \bigg)\,\bu(\x)-\bigg( \int_0^{\x}h \bu \bigg)\,\bv(\x).
\end{eqnarray*}

Every solution of \re{App3 1} is of the form $\,w=A\bu+B\bv+\bw$\, for some $(A,B)\in\R^2$. Since $h$\, is even, $\bw$\, is also even, so $w$ is even if and only if $A=0$. 

An even function $w=B\bv+\bw$\, is $2\pi$-periodic if and only if $w(\x+2\pi)-w(\x)=0$,\, that is, by \re{App3 diff},
\[
\bigg( \int_{\x}^{\x +2\pi}h \bv \bigg) \bu(\x)+\bigg[ \bigg( \int_0^{\x}h \bu \bigg) -B \bigg] \bar{V}^2 k \, \bu(\x) =0 \quad \forall \, \x .
\]

\noindent
We remove $\bu(\x)$,\, derive the expression with respect to $\x$\, and from \re{App3 diff} it results zero at any $\x$. Then the expression is a constant; we compute it at $\x=0$\, and obtain, since $h\bu$\, is odd and $2\pi$-periodic, that $w$\, is $2\pi$-periodic if and only if \,$B=\frac{1}{\bar{V}^2 k} \int_0^{2\pi}h \bv$. 

Thus, given $h$\, even $2\pi$-periodic, there exists a unique even $2\pi$-periodic  $w$\, such that $w''+\big( 3\b^2 +3\av{\b^2} \big) w =h$\, and this defines the operator $L$,
\begin{eqnarray} \label{App3 operatore L}
	L[h]=\bigg( \int_0^{\x}h \bv \bigg) \bu(\x)+ \bigg[
	\bigg( \frac{1}{\bar{V}^2 k}\int_0^{2 \pi}h \bv \bigg)\, - \, \int_0^{\x}h \bu \,\bigg] \,\bv(\x).
\end{eqnarray}

$L$\, is linear and continuous with respect to $\norm{\, \cdot \,}$; it is the Green operator of the equation \:$	x''+\big( 3\b^2 +3\av{\b^2} \big) x =h $,\;  so, by classical arguments, it is a bounded operator of \,${\mathcal H}_{\sigma,s}$\, into \,${\mathcal H}_{\sigma,s+2}$;\, the inclusion ${\mathcal H}_{\sigma,s+2} \hookrightarrow {\mathcal H}_{\sigma,s}$\, is compact, then $L:\hs \rightarrow \hs$\, is compact.  \str $\Box$

\jump 
\noindent \textbf{Lemma 3.} \emph{There holds the following equalities and inequalities:}

\jump
\noindent
\hspace{5pt}(L.1) \hspace{3pt}  $\av{\c^2}=\frac{1-2\m}{6\m}$\:  for $m=\m$. \: (Recall: $\c = \c (\, \cdot \, , m))$

\jump
\noindent
\hspace{5pt}(L.2) \hspace{3pt}  $\av{\frac{\s^2}{\d^2}}=\frac{1}{1-m}\av{\c^2}$\: for any $m$.

\jump
\noindent
\hspace{5pt}(L.3)  \hspace{3pt}  $m \av{\c^2 \frac{\s^2}{\d^2}}=1-2\av{\c^2}$\: for any $m$.

\jump
\noindent
\hspace{5pt}(L.4) \hspace{3pt} \emph{Exchange rule.}\ \   $\av{g L[h]}=\av{h L[g]} \ \ \,  \forall \, g,h$\, even $2\pi$-periodic.

\jump
\noindent
\hspace{5pt}(L.5) \hspace{3pt}  $1-3\av{\b^2 L[1]}=3\av{\b^2}\av{L[1]}$.

\jump
\noindent
\hspace{5pt}(L.6) \hspace{3pt} $\av{\b^2 L[\b]}=-\av{\b^2}\av{L[\b]}$.

\jump
\noindent
\hspace{5pt}(L.7) \hspace{3pt} $3\av{\b^2 L[\b^2]}= \av{\b^2} \,\Big( 1-3\av{L[\b^2]} \Big)$.

\jump
\noindent
\hspace{5pt}(L.8) \hspace{3pt} $\av{\b^2 L[\b]}=\av{\b L[\b^2]}=\av{L[\b]}=0$.

\jump
\noindent
\hspace{5pt}(L.9) \hspace{3pt}  $A_0 := \  1-3\av{\b^2 L[1]} \neq 0$.

\jump
\noindent
(L.10) \hspace{3pt}  $B_0 := \  1-6\av{\b L[\b]} \neq 0$ .

\jump
\noindent
(L.11) \hspace{3pt}  $C_0 := \ 1+6\av{\b L[\b]} \neq 0$ .

\jump
\noindent
(L.12) \hspace{3pt}  $A_0 \neq 1, \hspace{6pt} \av{L[\b^2]} \neq 0. $

\jump
\noindent
\emph{Proof.} 
(L.1) By construction of $\b$\, we have \,$\bar{\Om}^2(1-2\m)=3\av{\b^2}=3\bar{V}^2\av{\c^2(\cdot,\m)}$\: and \,$\bar{V}^2=2\m \bar{\Om}^2$,\, see Proof of Lemma 1.

\jump
\noindent (L.2) We observe that 
\[
\frac{d}{d\x}\left[ \frac{\c (\x)}{\d (\x)}\right]= \frac{(m-1) \s (\x)}{\d^2 (\x)},
\]

\noindent then we can integrate by parts
\[
\int_0^{4K}\frac{\s^2 (\x)}{\d^2 (\x)} \, d\x = \int_0^{4K} \frac{\s (\x)}{m-1}\, \frac{d}{d\x} \! \left[\frac{\c (\x)}{\d (\x)}\right]\, d\x=  \frac{1}{1-m} \int_0^{4K} \c^2 (\x) \, d\x.
\]

\jump
\noindent (L.3) We compute the derivative 
\[
\frac{d}{d\x} \left[ \frac{\c (\x) \s (\x)}{\d(\x)} \right]=2\c^2 (\x)-1+m\, \frac{\s^2 (\x) \c^2 (\x)}{\d^2 (\x)}
\]

\noindent and integrate on the period $[0,4K]$.

\jump
\noindent
(L.4) From the formula 
\re{App3 operatore L} of $L$\, we have 
	\[
	\begin{array}{rl}
\av{g L[h]}-\av{h L[g]}  
	= & \hspace{-4pt} 
	\av{\frac{d}{d\x} \big[ \big( \int_0^{\x}h \bv \big) \big( \int_0^{\x}g \bu \big) \big]} - \av{\frac{d}{d\x} \big[ \big( \int_0^{\x}h \bu \big) \big( \int_0^{\x}g \bv \big) \big]}+ \vspace{6pt} \\
 &  + \frac{1}{\bar{V}^2 k}\, 2\pi \big[ \av{h \bv}\av{g \bv} - \av{g \bv}\av{h \bv} \big]=0. 
   \end{array} 
\]

\jump
\noindent (L.5) By definition, $L[1]$\, satisfies \,$L[1]''+\big( 3\b^2+3\av{\b^2} \big) L[1]=1$,\, so we integrate on the period $[0,2\pi]$.

\jump
\noindent (L.6),(L.7) Similarly by definition of $L[\b],\ L[\b^2]$;\, recall that $\av{\b}=0$.

\jump
\noindent (L.8) By (L.6) and (L.4), it is sufficient to show that $\av{L[\b]}=0$.\, From the formula \re{App3 operatore L}, integrating by parts we have 
\[ 
\av{L[\b]}=-\av{\b \bv \Big( \int_0^\x \bu \Big)}-\av{\Big( \int_0^{\x} \b \bu \Big) \bv } + \frac{1}{\bar{V}^2 k} \Big( \int_0^{2\pi}\b \bv \Big) \av{\bv}.
\]

\noindent From the formulas \re{App3 2}, \re{App3 formula magica} of \,$\bu,\bv,$\,  recalling that \,$\b(\x)=\bar{V}\c(\bar{\Om} \x )$,\ we compute  
\begin{equation} \label{App4 1}   
\int_0^\x \bu \, = -\frac{1}{\bar{\Om}^2}\, \big( \c (\bar{\Om} \x)-1 \big), 
\qquad \int_0^\x \b  \bu \, = -\frac{\bar{V}}{2\bar{\Om}^2}\, \big( \c^2(\bar{\Om} \x) -1 \big).
\end{equation}

\noindent 
Observe that \,$\int_0^{2\pi}\c(\bar{\Om} \x) \frac{\s^2(\bar{\Om} \x)}{\d^2(\bar{\Om} \x)}\, d\x=0$\, by odd-symmetry with respect to $\frac{\pi}{2}$\, on $[0,\pi]$\, and periodicity. So, recalling that
$\bar{V}^2=2\m\bar{\Om}^2 $,\, we compute \,$\av{\bv}  =  \frac{\m k}{\pi}$.\, We can resume the computation of $\av{L[\b]}$\, obtaining
\[
\av{L[\b]}=\frac{3\bar{V}}{2 \bar{\Om}^2}\av{\bv(\x) \c^2(\bar{\Om} \x)} - \frac{\bar{V} \m k}{2 \pi \bar{\Om}^2}.
\]

\noindent
Since \,$\av{\c^3}=\av{\c^3 \frac{\s^2}{\d^2}}=0$\: by the same odd-symmetry reason, by \re{App3 formula magica} we have \,$\av{\bv(\x) \c^2(\bar{\Om} \x)}= \frac{\m k}{3\pi}$,\, and so  \,$\av{L[\b]}=0$. 

Moreover we can remark that by (L.4) there holds also $\av{\b L[1]}=0$.

\jump
\noindent (L.9) By (L.5), it is equivalent to show that $\av{L[1]} \neq 0$.\, From the formula \re{App3 operatore L}, integrating by parts we have
\[
 \av{L[1]}=\frac{2\pi}{\V^2 k} \av{\bv}^2-2\av{\Big( \int_0^\x \bu \Big) \bv }.
 \]
 
 \noindent We know that $\av{\bv}=\frac{\m k}{\pi}$,\: so by \re{App4 1}
 \[
  \av{L[1]}=\frac{1}{\Omb^2} \av{\bv(\x) \big( 2\c (\Omb \x) -1 \big) }.
 \]
 
\noindent 
>From the equalities (L.1) and (L.3) we have $\av{\bv(\x) \c (\Omb  \x)}=\frac{2}{3}(1-2\m)+\frac{\m k}{2\pi}$,\, thus 
\begin{equation} \label{App4 L[1]}
\av{L[1]}=\frac{4(1-2\m)}{3 \Omb^2}
\end{equation}

\noindent and this is strictly positive because $\m < \frac{1}{2}$.

\jump 
\noindent (L.10) From \re{App3 operatore L} integrating by parts we have
\begin{equation*} \label{App4 0,5}
\av{\b L[\b]}=-2\av{\b \bv \Big( \int_0^\x \b \bu \Big) }+\frac{2\pi}{\V^2 k} \av{\b \bv}^2.
\end{equation*}

\noindent Using (L.3), integrating by parts and recalling the definition \re{App3 def di k} of $k$\, we compute
\[
\av{\b \bv}= \V \m \av{\c^2}+\frac{\V \m k}{2\pi}+ \frac{\V (1-2\m)}{2}
\]

\noindent
and, by (L.1) and \re{App3 polinomio per k},
\begin{equation}  \label{App4 media beta v}
\av{\b \bv}= \frac{\V (7-8\m)}{12(1-\m)}.
\end{equation}

\noindent
By \re{App4 1}, $\av{\b \bv \Big( \int_0^\x \b \bu \Big) }=\frac{\V}{2 \Omb^2}\,\av{\b \bv \,\c^2} +\frac{\V}{2 \Omb^2}\,\av{\b \bv}$.\: The functions $\b$\, and $\bv$\, satisfy \,$\b''+\b^3+3\av{\b^2} \b=0$\, and $\bv''+3\b^2 \bv +3 \av{\b^2} \bv=0$,\, so that
\begin{equation} \label{App4 ecco l'errore}
\bv'' \b -\bv \b'' +2\b^3 \bv =0 .\vspace{4pt}
\end{equation}

\noindent
Deriving \re{App3 diff} we have \,$\bv'(2 \pi)-\bv'(0)=-\V^2 k$, so we can integrate \re{App4 ecco l'errore} obtaining
\[
\av{\b^3 \bv}=\frac{\V^3 k}{4\pi};
\]

\noindent
since $\av{\b^3 \bv }=\V^2 \av{\b \bv \,\c^2}$,\, we write
\[
\av{\b \bv \Big( \int_0^\x \b \bu \Big)}=-\frac{\m k}{4\pi}+\frac{\V}{2\Omb^2}\,\av{\b \bv}.
\]

\noindent Thus, by \re{App4 media beta v} and \re{App3 def di k}, we can express $\av{\b L[\b]}$\, in terms of $\m$\, only,
\begin{equation} \label{App4 beta L[beta] in m only}
\av{\b L[\b]}=\,\frac{32 \m^2 -32\m -1}{12(16\m^2 -16 \m+1)}=\:\frac{1}{6}-\frac{1}{4(16\m^2 -16 \m+1)}.
\end{equation}

\noindent 
The polynomial $p(m)=16m^2 -16 m+1$\: is non-zero for $m \in 
\big( \frac{2-\sqrt{3}}{4},\, \frac{2+\sqrt{3}}{4} \big)$\, and $\m \in (0.20, 0.21)$;\: so \,$B_0 =\frac{6}{4p(\m)}\neq 0$, in particular $B_0 \in (-1, -0.9)$.

\jump
\noindent (L.11) From \re{App4 beta L[beta] in m only} it follows that $C_0 \neq 0$,\, in particular \,$2.9 < C_0=2-\frac{3}{2p(\m)} < 3$.

\jump
\noindent (L.12) By Exchange rule (L.4) and (L.5), it is sufficient to show that $A_0 \neq 1$,\, that is \,$3\av{\b^2}\, \av{L[1]} \neq 1$.\, Recall that, by construction of $\m$,\, $3\av{\b^2}=\Omb^2 (1-2\m)$.\, So from \re{App4 L[1]} it follows
\[
 3\av{\b^2}\, \av{L[1]}=\frac{4}{3}\,(1-2\m)^2 ,
\]

\noindent and \,$\frac{4}{3}\,(1-2m)^2=1$\, if and only if \,$16m^2 -16 m+1=0$,\: while $\m \in (0.20, 0.21)$,\, like above; in particular \,$0.4\,<\, 3\av{\b^2}\, \av{L[1]} \,<\, 0.5.$     \str  $\Box$

\jump
\noindent
\textbf{\emph{Remark.}}  Approximated computations give 
\[
\begin{array}{lll}
\m \in (0.20,\, 0.21) \quad &   \bar{\sigma} \in (2.10, \,2.16)   \quad & \Omb \in (1.05,\, 1.06) \\
\V^2 \in (0.44, \,0.48)  \quad  & \av{\c^2} \in (2.85,\, 2.90)  \quad &  \av{\b^2} \in (1.27,\, 1.37).
\end{array}
\]

\jump
\noindent \textbf{Lemma 4.} \emph{The partial derivative \,$\partial_{Z}G(1,0,\b,\b)$\: is an invertible operator.}

\jump
\noindent
\emph{Proof.} Let $\,\partial_{Z}G(1,0,\b,\b)[\eta, h,k]=(0,0,0)\ $ for some $(\eta, h,k) \in Z,\,$ that is %\vspace{3pt}
\begin{eqnarray} \label{M 27}
  \begin{array}{l}
%\hspace{20pt} 
   6\eta \av{\b^2}+3\av{\b^2 h}+3\av{\b^2 k}=0 \vspace{6pt} \hspace{6pt}\\
% \hspace{20pt} 
  3\eta \big( \b^2 -\av{\b^2} \big) + h'' + \big( 3\b^2 +3\av{\b^2} \big) h   -3\av{\b^2 h} +6\av{\b k}\,\b =0 \vspace{6pt} \hspace{6pt}\\
% \hspace{20pt} 
  3\eta \big( \b^2 -\av{\b^2} \big) + k'' + \big( 3\b^2 +3\av{\b^2} \big) k -3\av{\b^2 k} +6\av{\b h}\,\b =0. \vspace{0pt} \hspace{6pt}
  \end{array}	
\end{eqnarray}

We evaluate the second and the third equation at the same variable and subtract; $\rho=h-k\,$ satisfies
\begin{eqnarray} \label{M 28}
	\rho''+\big( 3\b^2 +3\av{\b^2} \big)\rho -3\av{\b^2 \rho}-6\av{\b \rho}\,\b=0.
\end{eqnarray}
  
By definition of $L$,\, see Lemma 2, \re{M 28} can be written as
\begin{eqnarray} \label{M 29}
	\rho=3\av{\b^2 \rho}\,L[1] + 6\av{\b \rho}\,L[\b].
\end{eqnarray}

\noindent
Multiplying this equation by $\b^2 \,$ and integrating we obtain
\begin{eqnarray*} \label{M 30}
	\av{\b^2 \rho}\,\left( 1-3\av{\b^2 L[1]} \right) = 6 \av{\b \rho}\, \av{\b^2 L[\b]}.
\end{eqnarray*}

In Lemma 3\, we prove that $\big( 1-3\av{\b^2 L[1]} \big)=A_0\neq 0\,$ and $\av{\b^2 L[\b]}=0, \,$ then $\av{\b^2 \rho}=0$.

On the other hand, multiplying \re{M 29} by $\b\,$ and integrating we have
\begin{eqnarray*} \label{M 31}
  \av{\b \rho}\,\big( 1-6\av{\b L[\b]} \big) = 3 \av{\b^2 \rho}\, \av{\b L[1]};
\end{eqnarray*}

\noindent
in Lemma 3\, we show that $\big( 1-6\av{\b L[\b]} \big) =B_0\neq 0\,$ and 
$\av{\b L[1]}=0, \,$ then $\av{\b \rho}=0.\ $ From \re{M 29} we have so \,$\rho=0.\ $ Thus \,$h=k\,$ and \re{M 27} becomes 
\begin{eqnarray*} \label{M 32}
 \begin{array}{l}
\eta \av{\b^2}+\av{\b^2 h}=0 \vspace{6pt}\\
3\eta \big( \b^2 -\av{\b^2} \big) + h'' + \big( 3\b^2 +3\av{\b^2}     \big) h   -3\av{\b^2 h} +6\av{\b h}\,\b =0. 	
	 \end{array}
 \end{eqnarray*}
  
\noindent
By substitution we have 
\[
h=-3\eta L[\b^2]-6 \av{\b h}\,L[\b].
\]

Multiplying, as before, by $\b^2\,$ and by $\b \,$ and integrating, we obtain $\av{\b h}=\av{\b^2 h}=0 \,$ because $\big( 1+6\av{\b L[\b]} \big) =C_0\neq 0,$ \, $\av{\b L[\b^2]}=0,$ \, and $\, \av{\b^2}-3\av{\b^2 L[\b^2]}=3\av{\b^2}\, \av{L[\b^2]} \neq 0$, \ see Lemma 3\, again. Thus $h=0$, \, $\eta=0$ \, and the derivative  \,$\partial_{Z}G(1,0,\b,\b)\,$ is injective. 

The operator $Z \rightarrow Z,\: \: (\eta,h,k) \mapsto ((6\av{\b^2})^{-1} \eta,\,L[h],\,L[k])$\: is compact because $L$\, is compact, see Lemma 2.\, So, by the Fredholm Alternative, the partial derivative \,$\partial_{Z}G(1,0,\b,\b)$\, is also surjective.  \str $\Box$

\jump
By the Implicit Function Theorem and the regularity of $G$,\, using the rescaling \re{M 24 magic rescaling} we obtain, for $| b- \frac{1}{2} |$\, and \,$\e$\, small enough, the existence of a solution close to \,$(0,\b,\b)$\, for the $Z$-equation \re{M 13 projected problem}. 

More precisely: from Lemma 1\, and 4\, it follows the existence of a  $\mathcal{C}^1$-function $g$\, defined on a neighborhood of $\lambda=1$\, such that 
\[
G \big( \lambda, g(\lambda) \big)=0,
\]

\noindent
that is, $g(\lambda)$\, solves \re{M 25 il cuore},\, and \,$g(1)=(0,\b,\b)$.\, Moreover, for $|\lambda -1|$\, small, it holds
\begin{equation}  \label{stima g(lambda)}
\normm{g(\lambda) - g(1)} \leq \tilde{c} |\lambda -1|
\end{equation}

\noindent
for some positive constant $\tilde{c}$.\, In the following, we denote several positive constants with the same symbol $\tilde{c}$.
  
We set \,$\Phi_{(b,\e)} : (\uz,r,s) \mapsto (c,x,y)$\, the rescaling map \re{M 24 magic rescaling} and \,$H_{(b,\e)}:\R \times Z \rightarrow Z$\, the operator corresponding to the auxiliary bifurcation equation \re{M 22 aux Z eq}, which so can be written as 
\[
H_{(b,\e)}(\mu,z)=0. 
\]

We define 
\[
z^{*}_{(b,\e)}=\Phi_{(b,\e)}^{-1}[g(\lambda_{(b,\e)})],
\]

\noindent
thus it holds  \,$H_{(b,\e)}(0,z^{*}_{(b,\e)})=0$,\, that is, $z^{*}_{(b,\e)}$\, solves the bifurcation equation  \re{M 22 aux Z eq} for $\mu=0$.

We observe that $p_{(b,\e)}(0, z)=
\partial_z p_{(b,\e)}(0, z)=0$\, for every $z$\, and so, in particular, for $z=z_{(b,\e)}^*$;\, it follows that 
\begin{equation}  \label{M 33 derivata di H totale}
\partial_z  H_{(b,\e)}(0,z_{(b,\e)}^*)=(\Phi_{(b,\e)}^{-1})^3 \, \partial_z G\big( \lambda_{(b,\e)},\, g(\lambda_{(b,\e)})  \big) \,\Phi_{(b,\e)}.
\end{equation}

$G$\, is of class $\mathcal{C}^1$,\, so \,$\partial_Z G(\lambda,g(\lambda))$\, remains invertible for $\lambda$\, sufficiently close to 1.\, Notice that $\lambda_{(b,\e)}$\,  is sufficiently close to 1 if $| b- \frac{1}{2} |$\, and $\e$\, are small enough. Then, by \re{M 33 derivata di H totale}, the partial derivative \,$\partial_z  H_{(b,\e)}(0,z_{(b,\e)}^*)$\, is invertible. By the Implicit Function Theorem, it follows that for every $\mu$\, sufficiently small there exists a solution $z_{(b,\e)}(\mu)$\, of equation \re{M 22 aux Z eq}, that is
\[
H_{(b,\e)}(\mu, z_{(b,\e)}(\mu))=0. 
\]

We indicate $z_0=(0,\b,\b)$. The operators $\big( \partial_z H_{(b,\e)}(\mu,z) \big)^{-1}$\, and $\partial_\mu H_{(b,\e)}(\mu,z)$\, are bounded by some constant 
for every $(\mu,z)$\, in a neighborhood of $(0,z_0),$\, uniformly in $(b,\e)$,\, if $|b-1/2|,\e$\, are small enough.   
So the implicit functions $z_{(b,\e)}$\, are defined on some common interval $(-\mu_0,\mu_0)$\, for 
$|b-1/2|$,\, $\e$\, small, and it holds 
\begin{equation}  \label{1* pag.17}
\normm{z_{(b,\e)}(\mu)-z_{(b,\e)}^*} \leq \tilde{c} |\mu|	
\end{equation}

\noindent
for some $\tilde{c}$\, which does not depend on $(b,\e)$. 

Such a common interval $(-\mu_0,\mu_0)$\, permits the evaluation $z_{(b,\e)}(\mu)$\, at $\mu=\e$\, for $\e < \mu_0$,\, obtaining a solution of the original     bifurcation equation written in \re{M 13 projected problem}.  

Moreover,  \,$ \normm{ \Phi_{(b,\e)}^{-1} - \mathrm{Id}_Z } = |\sqrt{b(2+b\e^2)}-1| \leq |b-\frac{1}{2}| + \e^2$,\, so, by \re{stima g(lambda)} and triangular inequality, 
\begin{eqnarray}  \label{2* pag.17}
\begin{array}{c}
\normm{z_{(b,\e)}^* -z_0} \leq \tilde{c} ( |b-\frac{1}{2}|+\e^2 ). 
\end{array}
\end{eqnarray}

\noindent
Thus from \re{1* pag.17} and \re{2* pag.17} we have
\begin{eqnarray*}  \label{3* pag.17}
\begin{array}{c}
\normm{z_{(b,\e)}(\e) -z_0} \leq \tilde{c} ( |b-\frac{1}{2}|+\e ), 
\end{array}
\end{eqnarray*}

\noindent
and, by \re{bound for p}, 
\begin{eqnarray*}  \label{bound for p pag.17}
\begin{array}{c}
\normm{p \big( \e,\,z_{(b,\e)}(\e) \big)} \leq \tilde{c} \e.
\end{array}
\end{eqnarray*}

\noindent
\textbf{\emph{Remark.}} 
 Since the solutions $z_{(b,\e)}(\e)$\, are close to $z_0=(0,\b,\b)$,\, they actually depend on the two arguments $(\p,\pp)$;\, this is a necessary condition for the quasi-periodicity.

\jump
We define \,$
u_{(b,\e)}= z_{(b,\e)}(\e)+p_{(b,\e)}\big( \e,\,z_{(b,\e)}(\e) \big)$.\,  
Renaming $\mu_0 = \e_0$,\,  we have finally proved:

%Theorem 1

\jump
\noindent
\textbf{Theorem 1}. \emph{Let $\bar{\sigma} >0$,\, $\b$\, as in Lemma 1,\, $\tilde{B}_\gamma$\, as in \re{def di B_gamma nel testo} with $\gamma \in (0, \frac{1}{4})$.} 
\emph{For every $\sigma \in (0,\bar{\sigma})$,\,  
 there exist positive constants \,$\delta_0$,\, $\e_0$,\, $\bar{c}_1$,\, $\bar{c}_2$\, and the uncountable Cantor set}
\[
\bg=\Big\{ (b,\e) \in \Big(\frac{1}{2}- \delta_0,\, \frac{1}{2}+\delta_0 \Big) \times (0,\e_0 ) 
:  \ \frac{\e}{2+b\e^2},b\e^2 \in \tilde{B}_\gamma,\ \ 
 \frac{1+b\e^2}{\e} \notin \Q \Big\}
\]

\noindent
\emph{such that, for every $(b,\e) \in \bg $,\, there exists a solution $u_{(b,\e)} \in  \hs$\, of \re{M4}.} 
\emph{According to decomposition \re{M9 decomp2}, $u_{(b,\e)}$\, can be written as}
\[
  u_{(b,\e)}(\p,\pp)=\uz+r(\p)+s(\pp)+p(\p,\pp), 
\]

\noindent
\emph{where its components satisfy}  
\[
\begin{array}{c}
	\norm{r-\b}+ \norm{s-\b} +|\uz|   \leq \bar{c}_1 ( |b-\frac{1}{2}|+\e ), \qquad  \norm{p} \leq  \bar{c}_2 \e.
\end{array}
\]

\emph{As a consequence, problem \re{intro 1} admits uncountable many small amplitude, analytic, quasi-periodic solutions $v_{(b,\e)} $\, with two frequencies, of the form \re{intro 4}:}
\[ \begin{array}{rl}
    v_{(b,\e)}(t,x) \hspace{-6pt} & = \e\, u_{(b,\e)}\big( \e t, \: (1+b\e^2) t+x \big) \vspace{5pt}\\
      & = \e\,\big[\uz + r(\e t) +s((1+b\e^2) t+x ) + \mathcal{O}(\e) \big] \vspace{5pt}\\
      &  = \e\, \big[ \b (\e t) + \b ((1+b\e^2) t+x ) + \mathcal{O}\big(|b-\frac{1}{2}|+\e \big) \big]. 
\end{array}  
\]

\jump
\section{Waves traveling in opposite directions}

In this section we look for solutions of \re{intro 1} of the form \re{intro 2}, 
\[
v(t,x)=u(\om t+x, \, \omm t-x),
\]

\noindent
for \,$u \in \hs$.\, We introduce two parameters $(a,\e) \in \R^2$\, and set the frequencies as in \cite{Pro},
\begin{equation*} \label{L1}
\om =1+\e, \qquad \omm=1+a\e.
\end{equation*}

\noindent
For functions of the form \re{intro 2}, problem \re{intro 1} is written as
\[
L_{a,\e}[u] =  - u^3+ f (u)
\]

\noindent
where 
\[
L_{a,\e}=\e (2+\e) \,\partial_{\p}^2  + 2\big(2+(a+1)\e +a\e^2\big)\, \partial_{\p \pp}^2  + a\e (2+a\e)\, \partial_{\pp}^2 .
\]

We rescale \,$u \rightarrow \sqrt{\e}\, u$\, and define \,$f_{\e}(u)= \e^{-3/2}\,f(\sqrt{\e}\,u)$. Thus the problem can be written as
\begin{equation} \label{L3}
L_{a,\e}[u] =  -\e u^3+ \e f_\e (u). 
\end{equation}

For $\e=0$,\, the operator is \,$L_{a,0}=4 \partial_{\p \pp}^2$;\, its kernel is the direct sum \,$Z=C\oplus Q_1 \oplus Q_2$,\, see \re{M8 decomp}. 
Writing $u$\, in Fourier series we obtain an expression similar to \re{M5},
\begin{eqnarray*} \label{L4}
	L_{a,\e}[u]= \,- \hspace{-8pt} \sum_{(m,n)\in \Z ^2}D_{a,\e}(m,n)\, \hat{u}_{mn}\, e^{i m \p}\, e^{i n \pp},
\end{eqnarray*}
 
\noindent
where the eigenvalues \,$D_{a,\e}(m,n)$\, are given by
\begin{align*} %\label{L5}
  D_{a,\e}(m,n) \: = & \hspace{5pt} \e (2+\e)\, m^2 + a\e (2+a\e)\, n^2 + 2 \left( 2+(a+1)\,\e + a \e^2 \right) mn \nonumber \\
                = & \hspace{5pt} (2+\e)\,(2+a\e)\,\Big( m+\frac{a\e}{2+\e}\,n \Big) \Big( \frac{\e}{2+a\e}\, m +n \Big). 
\end{align*}                   

By Lyapunov-Schmidt reduction we project the equation \re{L3} on the four  subspaces,
\begin{eqnarray*} \label{L6 projected problem}
   \begin{array}{rl}
       0= & \hspace{-4pt} \uz ^3+3\uz\left(\av{r^2}+\av{s^2}\right)+\av{r^3}
                                                       +\av{s^3}+ \vspace{5pt}\\
          &+\Pi_C\left[(u^3-z^3)-f_{\e}(u)\right]\vspace{17pt}
                                \str  \left[ C-equation \right]\\
  -(2+\e)\,r''= & \hspace{-4pt}3\uz^2 r+3\uz \left(r^2 - \av{r^2}\right)+r^3 -\av{r^3}
                                                      +3\av{s^2}\,r+\vspace{5pt}\\            &     +\Pi_{Q_1}\left[(u^3-z^3)-f_{\e}(u)\right] \vspace{17pt}
                                  \str  \left[ Q_1-equation \right]\\
  -a(2+a\e)\, s'' = & \hspace{-4pt} 3\uz^2 s+3\uz \left( s^2 - \av{s^2} \right) 
                                        +s^3 -\av{s^3} + 3\av{r^2}\,s+\vspace{5pt}\\ 
          &    +\Pi_{Q_2}\left[(u^3-z^3)-f_{\e}(u)\right] \vspace{17pt}
                                 \str  \left[ Q_2 -equation \right]\\
  \lae [p]= & \hspace{-4pt} \e \, \Pi_P \left[-u^3 +f_{\e}(u)\right]. \vspace{5pt}
                           \str  \left[ P -equation \right] \vspace{5pt}\\
  \end{array}   
\end{eqnarray*}

We repeat the arguments of Appendix C\, and find a Cantor set $\ag$\, such that $|D_{a,\e}(m,n)|>\gamma$\, for every $(a,\e) \in \ag$.\, Then $L_{a,\e}$\, is invertible for $(a,\e) \in \ag$\, and the $P$-equation can be solved as in the section 4. 

We repeat the same procedure already shown in section 5 and solve the bifurcation equation. The only differences are: 

-	\,the parameter $a$\, tends to $1$\, instead of \,$b \rightarrow \frac{1}{2}$;

- \,the rescaling map is \,$\Psi_{(a,\e)}:(\uz,r,s) \mapsto (c,x,y)$,\, where		
\begin{eqnarray*} \label{L qualcosa}  
	\begin{array}{ll}
	r= \sqrt{2+\e}\ x  & \hspace{30pt} \uz=\sqrt{2+\e}\ c \\
	s= \sqrt{a(2+a\e)}\ y  & \hspace{30pt} \lambda=\lambda_{(a,\e)}=\frac{a(2+a\e)}{2+\e} , 
	\end{array}
\end{eqnarray*}

\hspace{9pt}instead of \,$\Phi_{(b,\e)}$\, defined in \re{M 24 magic rescaling}.  

We note that by means of the rescaling map \,$\Psi_{(a,\e)}$\, we obtain just the equation \re{M 25 il cuore}.
  Thus we conclude: 

\jump
\noindent
\textbf{Theorem 2}. \emph{Let $\bar{\sigma} >0$,\, $\b$\, as in Lemma 1,\, $\tilde{B}_\gamma$\, as in \re{def di B_gamma nel testo} with $\gamma \in (0, \frac{1}{4})$.} 
\emph{For every $\sigma \in (0,\bar{\sigma})$,\,  
 there exist positive constants \,$\delta_0$,\, $\e_0$,\, $\bar{c}_1$,\, $\bar{c}_2$\, and the uncountable Cantor set}
\[
\ag=\Big\{ (a,\e) \in (1- \delta_0,\, 1+\delta_0 ) \times (0,\e_0 ) 
:  \ \frac{a\e}{2+\e},\,\frac{\e}{2+a\e} \in \tilde{B}_\gamma,\ \ 
 \frac{1+\e}{1+a\e} \notin \Q \Big\} 
\]

\noindent
\emph{such that, for every $(a,\e) \in \ag $,\, there exists a solution $u_{(a,\e)} \in  \hs$\, of \re{L3}.} 
\emph{According to decomposition \re{M9 decomp2}, $u_{(a,\e)}$\, can be written as}
\[
  u_{(a,\e)}(\p,\pp)=\uz+r(\p)+s(\pp)+p(\p,\pp), 
\]

\noindent
\emph{where its components satisfy}  
\[
	\norm{r-\b}+ \norm{s-\b} +|\uz|   \leq \bar{c}_1 ( |a-1|+\e ), \qquad  \norm{p} \leq  \bar{c}_2 \e.
\]

\noindent
\emph{As a consequence, problem \re{intro 1} admits uncountable many small amplitude, analytic, quasi-periodic solutions $v_{(a,\e)}$ with two frequencies, of the form \re{intro 2}:}
\[
\begin{array}{rl} 
    v_{(a,\e)}(t,x) \hspace{-6pt} & = \sqrt{\e}\, u_{(a,\e)}\big( (1+\e) t+x, \: (1+a\e) t-x \big) \vspace{5pt}\\
    & = \sqrt{\e}\, \big[ \uz + r\big( (1+\e)t+x \big) + s\big( (1+a\e)t-x \big) + \mathcal{O}(\e) \big] \vspace{5pt} \\
      &  = \sqrt{\e}\, \big[ \b \big( (1+\e) t+x)\big) + \b \big( (1+a\e) t-x \big) + \mathcal{O}\big( |a-1|+\e \big) \big]. 
\end{array}
\]

\jump
\section{Appendix A. Hilbert algebra property of $\hs$} 
Let \,$u,v \in \hs$,\, $
u=\sum_{m \in \Z^2} \hat{u}_m \, e^{i m \cdot \ph}$,\:\: 
$v=\sum_{m \in \Z^2} \hat{v}_m \, e^{i m \cdot \ph}$.\,  
The product $uv$\, is
\[
uv= \sum_j \Big( \sum_k \hat{u}_{j-k} \hat{v}_k \Big)\, e^{i j \cdot \ph},
\]
 
\noindent
so its $\hs$-norm, if it converges, is 
\[
\norm{uv}^2= \sum_j \Big| \sum_k \hat{u}_{j-k} \hat{v}_k \Big|^2 \big(1+|j|^{2s}\big) \,e^{2|j| \sigma}.
\]
 
We define
\[
a_{jk}= \left[ \frac{ \big(1+|j-k|^{2s}\big)\, \big(1+|k|^{2s}\big) }{ \big(1+|j|^{2s}\big) } \right]^{1/2}.
\]

Given any $(x_k)_k$,\, it holds by H\"older inequality
\begin{equation}  \label{eq di x_k}
\Big| \sum_k x_k \Big|^2 = \Big| \sum_k \frac{1}{a_{jk}}\, x_k a_{jk} \Big|^2
 \leq c_j^2 \, \sum_k |x_k a_{jk} |^2 ,
\end{equation}

\noindent
where 
\[
c_j^2 := \sum_k \frac{1}{a_{jk}^2}. 
\]

We show that there exists a constant $c>0$\, such that \,$c_j^2 \leq c^2$\, for every \,$j \in \Z^2$.\, We recall that, fixed $p \geq 1$,\, it holds
\[
(a+b)^p \leq 2^{p-1}\, \big( a^p +b^p \big) \quad \forall \, a,b \geq 0.
\]

\noindent
Then, for $s \geq \frac{1}{2}$,\, we have
\begin{multline*}
1+|j|^{2s} \leq 1+ \big( |j-k| +|k| \big)^{2s} \leq 1+ 2^{2s-1}\, \big( |j-k|^{2s} + |k|^{2s} \big) \\
< 2^{2s-1}\, \big( 1+|j-k|^{2s} +1 +|k|^{2s} \big),
\end{multline*}

\noindent
so 
\[
\frac{1}{a_{jk}^2} < 2^{2s-1} \, \Big( \frac{1}{1+|j-k|^{2s} } + \frac{1}{1+|k|^{2s} } \Big).
\]  

\noindent
The series \,$\sum_{k \in \Z^2} \frac{1}{1+|k|^{p} }$\, converges if \,$p>2$,\, thus for $s>1$
\begin{multline*}
c_j^2 < 2^{2s-1}\, \Big( \sum_k \frac{1}{1+|j-k|^{2s} } + \sum_k \frac{1}{1+|k|^{2s} } \Big) = 2^{2s} \, \sum_{k \in \Z^2} \frac{1}{1+|k|^{2s} }\: := c^2 \,< \infty.
\end{multline*}

We put \,$x_k=\hat{u}_{j-k} \hat{v}_k   $\, in \re{eq di x_k} and compute
\begin{align*}
\Big| \sum_k  \hat{u}_{j-k} \hat{v}_k \Big|^2 \, \big( 1+|j|^{2s} \big) & \leq  
c^2 \, \sum_k |\hat{u}_{j-k} \hat{v}_k a_{jk} |^2 \, \big( 1+|j|^{2s} \big) \\
& = c^2 \, \sum_k |\hat{u}_{j-k} \hat{v}_k |^2 \, \big( 1+|j-k|^{2s} \big)\, \big( 1+|k|^{2s} \big), 
\end{align*}

\begin{align*}
\norm{uv}^2 & = \sum_j \Big| \sum_k \hat{u}_{j-k} \hat{v}_k \Big|^2 \big(1+|j|^{2s}\big) \,e^{2|j| \sigma}  \\
& \leq \sum_j c^2 \, \sum_k |\hat{u}_{j-k}|^2 | \hat{v}_k |^2 \, \big( 1+|j-k|^{2s} \big)\, \big( 1+|k|^{2s} \big)\, e^{2 ( |j-k| +|k| ) \sigma} \\
& = c^2 \, \sum_k \Big( \sum_j |\hat{u}_{j-k}|^2 \big( 1+|j-k|^{2s} \big)\,e^{2 |j-k| \sigma} \Big)\, | \hat{v}_k |^2  \big( 1+|k|^{2s} \big)\, e^{2 |k| \sigma} \\
& = c^2 \norm{u}^2 \norm{v}^2.
\end{align*}

So \,$\norm{uv}\leq c \norm{u}\norm{v}$\, for all $u,v \in \hs$.\, We notice that the constant $c$\, depends only on $s$,
\[
c= 2^s \bigg( \sum_{k \in \Z^2} \frac{1}{1+|k|^{2s}} \bigg)^{1/2}.
\]

\jump
\section{Appendix B. Change of form for quasi-periodic functions}

First an algebric proposition, then we show that one can pass from \re{intro 3} to \re{intro 4} without loss of generality.

\jump
\noindent
\textbf{Proposition.} \emph{Let $A,B \in \mathrm{Mat}_2 (\R)$ be invertible matrices such that $A B^{-1}$ has integer coefficient. Then, given any $u \in \hs$,\, the function  $v(t,x)=u\big(A(t,x) \big)$\,  can be written as \,$v(t,x)=w\big( B(t,x) \big)$\, for some \,$w \in \hs$, that is} \,$\{ u \cd A :\, u \in \hs \} \subseteq \{w \cd B : \, w \in \hs \}$.

\jump
\noindent
\emph{Proof.} Let $u \in \hs$.\, The function \,$u \circ A$\, belongs to \,$\{w \cd B : w \in \hs \}$\, if and only if \,$u \cd A \cd B^{-1} = w$\, for some \,$w \in \hs$,\, and this is true if and only if \,$u \cd AB^{-1}$\, is $2\pi$ periodic; since \,$AB^{-1} \in \mathrm{Mat}_2 (\Z)$,\, we can conclude. 
\str $\Box$

\jump
\noindent \textbf{Lemma.} \emph{The set of the quasi-periodic functions of the form  \re{intro 3} is equal to the set of the quasi-periodic functions of the form \re{intro 4}, that is,}
\begin{eqnarray*}
\begin{array}{l}
\Big\{ v : \,  v(t,x)=u(\om t+x,\,\omm t+x), \: 	(\om,\omm)\in \R^2,\,  
 \om \neq 0,\, \omm \neq 0,\, 
  \frac{\om}{\omm} \notin \Q,\,   u \in \hs \Big\} \str \vspace{6pt} \\
 = \Big\{ v : \,  v(t,x)=u\big( \e t, \,(1+b\e^2)t+x \big), \:  (b,\e) \in \R^2,\, \e \neq 0,\, (1+b\e^2) \neq 0, \str\\
  \str  \frac{1+b\e^2}{\e} \notin \Q,\,   u \in \hs \Big\}.
\end{array}
\end{eqnarray*}

\noindent 
\emph{Proof.} 
Given any $\om,\omm,b,\e$,\, we define 
\begin{eqnarray*} \label{App5 matrici A,B}
A= \bigg( \begin{array}{cc} 
  \om    & 1  \vspace{3pt} \\
  \omm    & 1 \\
   \end{array}  \bigg)
\quad \qquad 
B= \bigg( \begin{array}{cc} 
  \hspace{-4pt} \e    & 0  \vspace{3pt} \\
  \hspace{-4pt} ( 1+b\e^2)   & 1 \\
   \end{array}  \bigg).
\end{eqnarray*}   

Let $v(t,x)$ be any element of the set of quasi-periodic functions of the form  \re{intro 3}, that is \,$v=u \cd A$\, for some fixed  \,$\om,\omm \neq 0$\, such that $\frac{\om}{\omm} \notin \Q$\, and $u \in \hs$. We observe that $v$\, belongs to the set of quasi-periodic functions of the form \re{intro 4} if $v=w \cd B$\, for some $(b,\e)$\, such that $\e \neq 0,\: \frac{1+b\e^2}{\e} \notin \Q$\, and some $w \in \hs$.\, By the Proposition, this is true if we find $(b,\e)$\, such that $AB^{-1} \in  
\mathrm{Mat}_2 (\Z)$. We can choose
\[
b=\frac{\omm -1}{(\om -\omm )^2}, \qquad \e= \om -\omm,
\]

\noindent
so that\:$AB^{-1}= \Big( \begin{array}{cc} 
   1   & 1 \\
   0   & 1 \\
   \end{array}  \Big)$.
We notice that $\frac{1+b\e^2}{\e} \notin \Q$\, if and only if \,$\frac{\om}{\omm} \notin \Q$.

Conversely, we fix $(b,\e)$\, and look for $(\om,\omm)$\, such that $BA^{-1} \in \mathrm{Mat}_2 (\Z)$. This condition is satisfied if we choose the inverse transformation, \,$\om=1+\e+b\e^2$,\quad  $\omm= 1+b\e^2 $. 
\str $\Box$

\jump
\section{Appendix C. Small divisors}

Fixed $\gamma \in (0,\frac{1}{4})$,\, we have defined in \re{def di B_gamma nel testo} the set $\tilde{B}_\gamma\,$ of ``badly approximable numbers'' as
\begin{eqnarray*} \label{NA2 1 def di B_gamma}
	\tilde{B}_\gamma= \Big\{ x\in \R : \left| m+nx \right| > \frac{\gamma}{| n | } \ \ \forall \, m,n \in \Z,\, m\neq 0,\,n \neq 0 \Big\}.
\end{eqnarray*}

$\tilde{B}_\gamma\,$ is non-empty, symmetric, it has zero Lebesgue-measure and it  accumulates to 0. 
Moreover, for every $\delta >0$,\, $\tilde{B}_\gamma \cap(-\delta, \delta)$ is uncountable. 

In fact, $\tilde{B}_\gamma\,$ contains all the irrational numbers whose continued fractions expansion is of the form $[0,a_1, a_2,  \dots ]$,\, with \,$a_j < \gamma^{-1}-2$\: for every  $j \geq 2$. Such a set is uncountable: since $\gamma^{-1} -2 >2$,\, for every $j\geq 1$\, there are at least two choices for the value of $a_j$. Moreover, it accumulates to 0: if $y=[0,a_1,a_2, \dots]$, it holds $0<y<a_1^{-1}$,\, and $a_1$\, has no upper bound. See also Remark 2.4 in \cite{BP} and, for the inclusion of such a set in $\tilde{B}_\gamma$,\, the proof of Theorem 5F in \cite{Schmidt}, p. 22.

We prove the following estimate.

\jump
\noindent  
\textbf{Proposition.} \emph{Let $\gamma \in (0,\frac{1}{4}), \, \delta \in (0,\frac{1}{2}). \,$ Then for all $x, y \in \tilde{B}_\gamma\cap(-\delta, \delta)\,$ it holds      }
\begin{eqnarray*} \label{NA2 2 prop}
	\big| (m+nx)(my+n) \big| >\gamma (1-\delta -\delta^2) \quad \forall \, m,n \in \Z,\, m,n \neq 0. \vspace{7pt}
\end{eqnarray*}

\noindent \emph{Proof.} We shortly set $D=\left| (m+nx)(my+n)\right|.\,$ There are four cases.

\noindent Case 1. $|m+nx|>1, \ |my+n|>1. \  $ Then $| D |>1$. 

\noindent Case 2. $|m+nx|<1, \ |my+n|>1. \  $ 
 Multiplying the first inequality by $|y|,$
\begin{eqnarray*}
  \begin{array}{rl}
|y|> & \hspace{-4pt} |my+nxy|\,= \, |my+n-n(1-xy)| \vspace{3pt} \\
\geq & \hspace{-4pt} \big| |n(1-xy)| - |my+n| \big| \,\geq  \, |n(1-xy)| - |my+n|, \vspace{3pt}
  \end{array}
\end{eqnarray*}

\noindent so \,$|my+n|>|n|(1-xy) - |y|$\, and 
\begin{eqnarray*}
  \begin{array}{rl}
  |D|> & \hspace{-4pt} \frac{\gamma}{|n|} \big[ |n|(1-xy)-|y| \big] = 
  \gamma \Big[ (1-xy) - \frac{|y|}{|n|} \Big] \vspace{3pt} \\
  > & \hspace{-4pt} \gamma [(1-\delta^2)-\delta].     \vspace{3pt} 
  \end{array}
\end{eqnarray*}

\noindent Case 3. $|m+nx|>1, \ |my+n|<1. \  $ Analogous to case 2.

\noindent Case 4. $|m+nx|<1, \ |my+n|<1. \  $ Dividing the first inequality by $|n|,\,$ for triangular inequality we have 
\[
\left| \frac{m}{n} \right| \,\leq\, \left| \frac{m}{n}+x \right| + |x| \,<\,   \frac{1}{|n|}+\delta ,
\]

\noindent and similarly $\left| \frac{n}{m} \right| < \frac{1}{|m|}+\delta$. So
\[
\Big( \frac{1}{|n|}+\delta \Big) \Big( \frac{1}{|m|}+\delta \Big) > 
\left|  \frac{n}{m} \cdot \frac{m}{n} \right| =1.
\]

\noindent If \,$|n|, \, |m| \geq 2, \,$ then \,$\left( \frac{1}{|n|}+\delta \right) \left( \frac{1}{|m|}+\delta \right)<1,\,$ a contradiction. It follows that at least one between $|n|$ and $|m|$ is equal to 1. Suppose $|n|=1$.\, Then \,$|m+nx|=|m\pm x|\geq |m|-\delta \,$ and
\[
|D|>\frac{\gamma}{|m|} \big( |m|-\delta \big) = \gamma \Big( 1-\frac{\delta}{|m|} \Big) \geq \gamma (1-\delta).
\]
 
\noindent If $|m|=1\,$ the conclusion is the same. \str $\Box$

\jump
Fixed $\gamma \in (0,\frac{1}{4})$\, and $\delta \in (0,\frac{1}{2})$,\, we define the set 
\begin{eqnarray*} \label{NA2 3 def di Bgdelta}
  \begin{array}{rl}
	B(\gamma,\delta)=\Big\{ (b,\e) \in \R^2 : & \e \neq 0, \ \   1+b\e^2 \neq 0, \ \    
	2+b\e^2 \neq 0,    \vspace{2pt} \\
	& \frac{1+b\e^2}{\e} \notin \Q, \ \  \frac{\e}{2+b\e^2},\, b\e^2  \in \tilde{B}_\gamma \cap (-\delta, \delta) \Big\}
  \end{array}
\end{eqnarray*}

\noindent
and the map 
\[
g: B(\gamma,\delta) \rightarrow \R^2, \qquad 
 g(b,\e)=\Big( \frac{\e}{2+b\e^2}, \, b\e^2 \Big).
\]
 
\noindent  $g'(b,\e)\,$ is invertible on $B(\gamma,\delta).\,$ Its imagine is the set \,$R(g)=\{(x,y)\in \R^2 : x,y \in \tilde{B}_\gamma\cap (-\delta, \delta),\ \frac{1}{x}-y \notin \Q \}$\, and its inverse is
\begin{eqnarray*}
	g^{-1}(x,y)= \left( \frac{y(1-xy)}{2x},\  \frac{2x}{1-xy} \right).
\end{eqnarray*}
  
Thus $B(\gamma,\delta)\,$ is homeomorphic to \,$R(g)=\{(x,y)\in \tilde{B}_\gamma^2 : |x|,|y| <\delta,\ \frac{1}{x}-y \notin \Q \}.$\, We observe that, fixed any $\bar{x} \in \tilde{B}_\gamma \cap (-\delta, \delta)$, it occurs $\frac{1}{\bar x }-y \in \Q$\, only for countably many numbers $y$.\, We know that $\tilde{B}_\gamma \cap (-\delta, \delta)$\, is uncountable so, removing from  \,$[\tilde{B}_\gamma \cap (-\delta, \delta)]^2$\,     the couples $\{ (\bar x ,y) : y=\frac{1}{\bar x }-q\ \, \exists \,q\in \Q \}$,\, it remains uncountably many other couples. Thus $R(g)$\, is uncountable and so, through $g$,\, also $B(\gamma,\delta)$.\

Moreover, if we consider couples $(x,y) \in [\tilde{B}_\gamma \cap (-\delta, \delta)]^2\,$ such that $x \rightarrow 0\,$ and $(x/y) \rightarrow 1,\,$ applying $g^{-1}\,$ we find couples $(b,\e) \in B(\gamma,\delta)\,$ which satisfy  \,$\e\rightarrow 0$,\: $b\rightarrow 1/2$. In other words, the set \,$B(\gamma,\delta)$\, accumulates to $(1/2,0)$.

Finally we estimate $D_{b,\e}(m,n)\,$ for $(b,\e) \in B(\gamma,\delta)$. We have 
\[
|2+b\e^2|=  \frac{2}{|1-xy|}  > \frac{2 }{1+\delta^2},
\]
 \noindent 
so from the previous Proposition and \re{M6} it follows 
\[
\big| D_{b,\e}(m,n) \big| = |D|\,|2+b\e^2| > \gamma (1-\delta -\delta^2)  \frac{2}{1+\delta^2} .
\] 

\noindent
The factor on the right of $\gamma \,$ is greater than 1 if we choose, for example, \,$\delta=1/4; \,$ we define \,$\bg = B(\gamma,\delta)_{|\delta=\frac{1}{4}}$,\, so that there holds
\[
\big| D_{b,\e}(m,n) \big|  > \gamma \qquad \forall \, (b,\e) \in \bg.
\]  

We can observe that the condition \,$\frac{1+b\e^2}{\e} \notin \Q$\, implies \,$1+b\e^2 \neq0$,\, that $\frac{\e}{2+b\e^2} \in \tilde{B}_\gamma $\, implies \,$\e \neq 0$\, and \,$|b\e^2|<\delta$\, implies \,$2+b\e^2 \neq0$,\, so that we can write
\begin{equation*}  \label{NA2 vero bg}
\bg=\Big\{ (b,\e) \in \R^2 : \ \frac{\e}{2+b\e^2},b\e^2 \in \tilde{B}_\gamma,\ \ 
\Big| \frac{\e}{2+b\e^2} \Big|, |b\e^2 | < \frac{1}{4}, \ \ \frac{1+b\e^2}{\e} \notin \Q \Big\} . 
\end{equation*} 

We notice also that, for \,$| b- \frac{1}{2} |$\, and  \,$\e$\: small enough, there holds  automatically \,$| \frac{\e}{2+b\e^2} | < \frac{1}{4}$,\, $|b\e^2 | < \frac{1}{4}$.\: So, if we are interested to couples $(b,\e)$\, close to $\big( \frac{1}{2},0 \big)$,\, say $|b-\frac{1}{2}|< \delta_0$,\, $|\e|< \e_0$,\, we can write 
\begin{equation*}  \label{NA2 vero bg ristretto}
\bg=\Big\{ (b,\e) \in \Big(\frac{1}{2}- \delta_0,\, \frac{1}{2}+\delta_0 \Big) \times (0,\e_0 ) 
: \ \frac{\e}{2+b\e^2},b\e^2 \in \tilde{B}_\gamma,\ \ 
 \frac{1+b\e^2}{\e} \notin \Q \Big\} .
\end{equation*}

%%NUMERACCI
%\begin{eqnarray*} \begin{array}{ll}
%\Om \in (1.056534,\:  1.059845) & 
%\quad V^2 \in (0.446506,\:  0.471774) \vspace{4pt}\\
%\av{\c^2} \in (2.857143,\:  2.90) &
%\quad \av{\b^2} \in (1.275731,\:  1.368145) \\ 
%\end{array}
%\end{eqnarray*}
%
%\noindent 
%\dots e poi li ho approssimati per scriverli nel Remark.
%

\jump

\jump

\end{document}